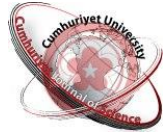

# Notes on Prime Near-Rings with Multiplicative Derivation

Zeliha BEDIR[1*], Oznur GOLBASI

[1]*Cumhuriyet University, Faculty of Science, Department of Mathematics, 58140 Sivas,*

**Abstract:** Let $N$ be a left near ring. A map $d: N \to N$ is called a nonzero multiplicative derivation if $d(xy) = xd(y) + d(x)y$ holds for all $x, y \in N$. In the present paper, we shall extend some well known results concerning commutativity of prime rings for nonzero multiplicative derivations of a left prime near-ring $N$.

## Çarpımsal Türevli Asal Yakın Halkalar Üzerine Notlar

**Özet**: $N$ bir sol yakın halka olsun. $d: N \to N$ dönüşümü her $x, y \in N$ için $d(xy) = xd(y) + d(x)y$ koşulunu sağlıyorsa $d$ ye bir çarpımsal türev denir. Bu makalede, asal halkalarda iyi bilinen bazı komütatiflik koşulları, çarpımsal türevli sol asal yakın halkalar için genelleştirilecektir.

**Anahtar Kelimeler:** Asal halka, yakın halka, türev, çarpımsal türev

## 1. INTRODUCTION

An additively written group $(N, +)$ equipped with a binary operation $. : N \to N, (x, y) \to xy$, such that $x(yz) = (xy)z$ and $x(y+z) = xy + xz$ for all $x, y, z \in N$ is called a left near-ring. A near-ring $N$ is called zero symmetric if $0x = 0$ for all $x \in N$ (recall that left distributive yields $x0 = 0$). A near-ring $N$ is said to be 3-prime if $xNy = \{0\}$ implies $x = 0$ or $y = 0$. For any $x, y \in N$, as usual $[x, y] = xy - yx$ and $xoy = xy + yx$ will denote the well-known Lie and Jordan products respectively. The set $Z = \{x \in N \mid yx = xy \text{ for all } y \in N\}$ is called multiplicative center of $N$. A mapping $d: N \to N$ is said to be a derivation if $d(xy) = xd(y) + d(x)y$ for all $x, y \in N$. $N$ is said to be $2$-torsion free if $x \in N$ and $x + x = 0$ implies $x = 0$.

Since Posner published his paper [11] in 1957, many authors have investigated properties of derivations of prime and semiprime rings. The study of derivations of near-rings was initiated by Bell and Mason in 1987 [1]. There has been a great deal of work concerning commutativity of prime and semiprime rings and near-rings with derivations satisfying with certain differential identities. (see references for a partial bibliography).

___

* Corresponding author. *Email address:* zelihabedir@cumhuriyet.edu.tr

# Notes on Prime Near-Rings with Multiplicative Derivation

In [7], Herstein has proved that if $R$ is a prime ring of characteristic different from 2 and if $d$ is a nonzero derivation of $R$ such that $d(R) \subseteq Z$, then $R$ is commutative. In [3], Bell and Kappe have proved that $d$ is a derivation of $R$ which is either a homomorphism or an anti-homomorphism in semiprime ring $R$ or a nonzero right ideal of $R$ then $d = 0$. In [5], Daif and Bell proved that if $R$ is semiprime ring, $U$ is a nonzero ideal of $R$ and $d$ is a derivation of $R$ such that $d([x, y]) = \pm[x, y]$, for all $x, y \in U$, then $U \subseteq Z$. All of these results were extended to near rings.

In [4], the notion of multiplicative derivation was introduced by Daif motivated by Martindale in [8]. $d : R \to R$ is called a multiplicative derivation if $d(xy) = xd(y) + d(x)y$ holds for all $x, y \in R$. These maps are not additive. In [6], Goldman and Semrl gave the complete description of these maps. We have $R = C[0,1]$, the ring of all continuous (real or complex valued) functions and define a map $d : R \to R$ such as

$$d(f)(x) = \begin{cases} f(x) \log |f(x)|, & f(x) \neq 0 \\ 0, & \text{otherwise} \end{cases}.$$

It is clear that $d$ is multiplicative derivation, but $d$ is not additive.

Recently, some results concerning commutativity of prime rings with derivations were proved for multiplicative derivations. It is natural to look for comparable results with multiplicative derivations of near-rings. In the present paper, we shall extend above mentioned results for multiplicative derivations of 3-prime near-ring $N$. Also, we will prove some commutativity conditions.

**Chapter 1:**

**Lemma 1** *[2, Lemma 1.2] Let $N$ be a $3-$prime near-ring.*

$(i)$ *If $z \in Z \setminus \{0\}$, then $z$ is not a zero divisor.*

$(ii)$ *If $Z$ contains a nonzero element $z$ for which $z + z \in Z$, then $(N, +)$ is abelian.*

$(iii)$ *If $z \in Z \setminus \{0\}$ and $x \in N$ such that $xz \in Z$ or $zx \in Z$, then $x \in Z$.*

**Lemma 2** *[2, Lemma 1.5] Let $N$ be a $3-$prime near ring. If $Z$ contains a nonzero semigroup ideal of $N$, then $N$ is commutative ring.*

**Lemma 3** *[9, Lemma 2.1] A near-ring $N$ admits a multiplicative derivation if and only if it is zero symmetric.*

**Lemma 4** *Let $N$ be a near-ring and $d : N \to N$ multiplicative derivation of $N$. Then*

$(xd(y) + d(x)y)z = xd(y)z + d(x)yz,$ for all $x, y, z \in N.$

*Proof:* By calculating $d(xyz)$ in two different ways, we see that

$$d((xy)z) = xyd(z) + d(xy)z$$

and

$$d(x(yz)) = xd(yz) + d(x)yz$$
$$= xyd(z) + xd(y)z + d(x)yz.$$

Hence we have

$$d(xy)z = xyd(z) + xd(y)z$$

and so

$$(xd(y) + d(x)y)z = xd(y)z + d(x)yz, \text{for all } x, y, z \in N.$$

**Lemma 5** *Let $N$ be a 3-prime near-ring and $a \in N$. If $N$ admits a nonzero multiplicative derivation $d$ such that $d(N)a = 0$ (or $ad(N) = 0$), then $a = 0$.*

*Proof.* By the hypothesis, we get

$$d(xy)a = 0, \text{for all } x, y \in N.$$

Expanding this equation with Lemma 4 and using the hypothesis, we have

$$d(x)Na = (0), \text{for all } x \in N.$$

Since $N$ is 3-prime near-ring and $d \neq 0$, we obtain that $a = 0$.

$ad(N) = 0$ can be proved by applying the same techniques.

**Theorem 1** *Let $N$ be a 3-prime near-ring. If $N$ admits a nonzero multiplicative derivation $d$ such that $d(N) \subseteq Z$, then $N$ is a commutative ring.*

*Proof.* For any $x, y \in N$, we get $d(xy) \in Z$, and so

$$d(xy)y = yd(xy).$$

That is

$$(xd(y) + d(x)y)y = y(xd(y) + d(x)y).$$

Using Lemma 4, we get

$$xd(y)y + d(x)yy = yxd(y) + yd(x)y.$$

Since $d(N) \subseteq Z$, we arrive at

$$d(y)xy + d(x)yy = d(y)yx + d(x)yy$$

and so

$$d(y)[x, y] = 0.$$

Using Lemma 1 (i), we have for each fixed $y \in N$ either $d(y) = 0$ or $y \in Z$.

Now, we assume $d(y) = 0$. For any $x \in N$, we have $d(xy) \in Z$ by the hypothesis. Since $d(y) = 0$, we get $d(xy) = d(x)y \in Z$, for all $x \in N$. By Lemma 1 (iii), we get $d(x) = 0$, for all $x \in N$ or $y \in Z$. Since $d \neq 0$, we must have $y \in Z$. Hence we arrive at $y \in Z$ for any cases. That is $N \subseteq Z$, and so $N$ is commutative near-ring by Lemma 2.

**Theorem 2** *Let $N$ be a 3-prime near-ring and $d$ a multiplicative derivation of $N$ such that $d(xy) = d(x)d(y)$, for all $x, y \in N$, then $d = 0$.*

*Proof.* In view of our hypothesis, we have

$$xd(y) + d(x)y = d(x)d(y), \text{ for all } x, y \in N. \tag{2.1}$$

Replacing $y$ by $yz$ in (2.1), we get

$$xd(yz) + d(x)yz = d(x)d(yz).$$

By our hypothesis, we have

$$xd(y)d(z) + d(x)yz = d(x)d(y)d(z)$$

and so

$$xd(y)d(z) + d(x)yz = d(xy)d(z).$$

Since $d$ is multiplicative derivation of $N$, we arrive at

$$xd(y)d(z) + d(x)yz = (xd(y) + d(x)y)d(z).$$

By Lemma 4, we get

$$xd(y)d(z) + d(x)yz = xd(y)d(z) + d(x)yd(z), \text{ for all } x, y, z \in N.$$

That is

$$d(x)yz = d(x)yd(z), \text{for all } x, y, z \in N.$$

Since $N$ is left near-ring, we have

$$d(x)N(d(z) - z) = (0), \text{for all } x, z \in N.$$

By the 3-primeness of $N$, we arrive at

$$d = 0 \text{ or } d(z) = z, \text{for all } z \in N.$$

If $d(z) = z$, for all $z \in N$, then

$$d(xy) = xd(y) + d(x)y$$

$$xy = xy + xy$$

$$xy = 0, \text{for all } x, y \in N.$$

This yields that $N = (0)$, a contradiction. So, we must have $d = 0$. This completes the proof of our theorem.

**Theorem 3** *Let $N$ be a 3-prime near-ring and $d$ a multiplicative derivation of $N$ such that $d(xy) = d(y)d(x)$, for all $x, y \in N$, then $d = 0$.*

*Proof.* By our hypothesis, we have

$$xd(y) + d(x)y = d(y)d(x), \text{for all } x, y \in N. \qquad (2.2)$$

Replacing $y$ by $xy$ in (2.2), we get

$$xd(xy) + d(x)xy = d(xy)d(x).$$

In view of our hypothesis, we have

$$xd(y)d(x) + d(x)xy = d(xy)d(x).$$

Using $d$ is multiplicative derivation of $N$, we arrive at

$$xd(y)d(x) + d(x)xy = (xd(y) + d(x)y)d(x).$$

By Lemma 4, we get

$$xd(y)d(x) + d(x)xy = xd(y)d(x) + d(x)yd(x)$$

and so

$$d(x)xy = d(x)yd(x), \text{ for all } x, y \in N. \tag{2.3}$$

Taking $yz$ instead of $y$ in (2.3) and using (2.3), we obtain that

$$d(x)N[z, d(x)] = 0, \text{ for all } x, z \in N.$$

By the 3-primeness of $N$, we get

$$d(x) = 0 \text{ or } d(x) \in Z.$$

Now, $d(x) = 0$ implies that $d(x) \in Z$. So, we have $d(N) \subseteq Z$ for any cases. By Theorem 1, we obtain that $N$ is commutative ring or $d = 0$. If $N$ is commutative ring, then $d(xy) = d(y)d(x) = d(x)d(y)$, for all $x, y \in N$. Hence, we get $d = 0$ by Theorem 2. This completes the proof.

**Theorem 4** *Let $N$ be a 3-prime near-ring and $d$ a nonzero multiplicative derivation of $N$ such that $d([x, y]) = [d(x), y]$, for all $x, y \in N$, then $N$ is commutative ring.*

*Proof.* Replacing $xy$ instead of $y$ in the hypothesis, we get

$$d(x[x, y]) = [d(x), xy].$$

Expanding this equation and using the hypothesis, we have

$$xd([x, y]) + d(x)[x, y] = [d(x), xy]$$

$$x[d(x), y] + d(x)[x, y] = [d(x), xy]$$

$$xd(x)y - xyd(x) + d(x)[x, y] = d(x)xy - xyd(x).$$

On the other hand, replacing $y = 0$ in the hypothesis, we arrive at $d(0) = 0$. Again replacing $x$ instead of $y$ in the hypothesis, we get

$$[d(x), x] = 0$$

and so

$$d(x)x = xd(x), \text{ for all } x \in N.$$

Now, using this in the above equation, we find that

$$d(x)xy - xyd(x) + d(x)[x, y] = d(x)xy - xyd(x)$$

$$d(x)[x, y] = 0$$

and so

$$d(x)xy = d(x)yx, \text{ for all } x, y \in N.$$

Replacing $y$ by $yz$ in this equation and using this, we have

$$d(x)N[x, z] = (0), \text{ for all } x, z \in N.$$

This yields that

$$d(x) = 0 \text{ or } x \in Z.$$

If $d(x) = 0$, then $d(x) \in Z$. On the otherwise, if $x \in Z$ then $[d(x), y] = 0$, for all $y \in N$ by the hypothesis. Hence we have $d(x) \in Z$. Thus we arrive at $d(x) \in Z$, for both cases. That is $d(N) \subseteq Z$, and so, we obtain that $N$ is commutative ring by Theorem 1.

**Theorem 5** *Let $N$ be a 3-prime near-ring and $d$ a nonzero multiplicative derivation of $N$ such that $[d(x), y] = [d(x), d(y)]$, for all $x, y \in N$, then $N$ is commutative ring.*

*Proof.* If $d(x) \in Z$, then there is nothing to prove. So we assume that $d(x) \notin Z$, for any $x \in N$. In the view of the hypothesis, we get

$$[d(x), y] = [d(x), d(y)], \text{ for all } x, y \in N.$$

Writing $d(x)y$ instead of $y$ in this equation, we get

$$[d(x), d(d(x)y)] = [d(x), d(x)y]$$

$$d(x)d(d(x)y) - d(d(x)y)d(x) = d(x)[d(x), y].$$

Using $d$ is multiplicative derivation of $N$ and Lemma 4, we arrive at

$$d(x)d(x)d(y) + d(x)d^2(x)y - (d(x)d(y)d(x) + d^2(x)yd(x)) = d(x)[d(x), y].$$

By the hypothesis, we have

$$d(x)d(x)d(y) + d(x)d^2(x)y - (d(x)d(y)d(x) + d^2(x)yd(x)) = d(x)[d(x), d(y)].$$

Expanding this term and using $-(a+b) = -b-a,$ we arrive at

$$d(x)d(x)d(y) + d(x)d^2(x)y - d^2(x)yd(x) - d(x)d(y)d(x) = d(x)d(x)d(y) - d(x)d(y)d(x)$$

and so

$$d(x)d^2(x)y = d^2(x)yd(x), \text{ for all } x, y \in N.$$

Replacing $yz$ instead of $y$ in the last equation, we find that

$$d^2(x)N[d(x), z] = (0), \text{ for all } x, z \in N.$$

By the 3-primeness of $N$, we get for each $x \in N$

$$d^2(x) = 0 \text{ or } d(x) \in Z.$$

Since $d(x) \notin Z,$ we must have $d^2(x) = 0,$ for all $x \in N.$ Writing $d(y)$ instead of $y$ in the hypothesis and using $d^2(y) = 0,$ we arrive at $[d(x), d(y)] = 0.$ Again using this in the hypothesis, we have $[d(x), y] = 0,$ and so $d(x) \in Z,$ a contradiction. Hence, we must have $d(N) \subseteq Z,$ and so, $N$ is commutative ring by Theorem 1. This completes the proof.

**Theorem 6** *Let $N$ be a 3-prime near-ring, $d$ a multiplicative derivation of $N$. If $[d(x), y] \in Z,$ for all $x, y \in N,$ then $N$ is a commutative ring.*

*Proof.* Replacing $y$ by $d(x)y$ in the hypothesis yields that

$$d(x)[d(x), y] \in Z, \text{ for all } x, y \in N.$$

By Lemma 1 (iii), we get

$$d(x) \in Z \text{ or } [d(x), y] = 0, \text{ for all } x, y \in N.$$

For any cases, we obtain that $d(N) \subseteq Z.$ By Theorem 1, we obtain that $N$ is a commutative ring.

**Acknowledgments:** This work is supported by the Scientific Research Project Fund of Cumhuriyet University under the project number F-496.

**REFERENCES**
[1]. Bell, H. and Mason, G., On derivations in near rings, *Near rings and Near fields, North-Holland Mathematical Studies,* 137, 31-35, (1987).
[2]. Bell, H. E., On derivations in near-rings II, *Kluwer Academic Pub. Math. Appl., Dordr.*, 426, 191-197, (1997).

[3]. Bell, H. E. and Kappe, L. C., Rings in which derivations satisfy certain algebraic conditions, *Acta Math. Hungarica,* 53, 339-346, (1989).

[4]. Daif, M. N., When is a multiplicative derivation additive, *Int. J. Math. Math. Sci.,* 14(3), 615-618, (1991).

[5]. Daif, M. N. and Bell, H. E., Remarks on derivations on semiprime rings, *Int. J. Math. Math. Sci.,* 15(1), 205-206, (1992).

[6]. Goldman, H. and Semrl, P., Multiplicative derivations on $C(X)$, *Monatsh Math.*, 121(3), 189-197, (1969).

[7]. Herstein, I. N., A note on derivations, *Canad. Math. Bull.,* 21(3), 369-370, (1978).

[8]. Martindale III, W. S., When are multiplicative maps additive, *Proc. Amer. Math. Soc.*, 21, 695-698, (1969).

[9]. Kamal, A. M. and Al-Shaalan, K. H., Existence of derivations on near-rings, *Math. Slovaca,* 63, no:3, 431-438, (2013).

[10]. Koç, E. and Gölbaşı, Ö., Semiprime near-rings with multiplicative generalized $(\theta,\theta)-$derivations, *Fasciculi Mathematici,* 57, 105-119, (2016).

[11]. Posner, E. C. Derivations in Prime Rings, *Proc. Amer. Math. Soc.,* 8, 1093-1100, (1957).